\documentclass{article}
\usepackage{amssymb}
\usepackage{amsmath}

\newtheorem{theorem}{Theorem}

\newenvironment{proof}[1][Proof]{\noindent\textbf{#1.} }{\ \rule{0.5em}{0.5em}}
\input{tcilatex}
\begin{document}

\title{The Jacobsthal and Jacobsthal-Lucas sequences associated with pseudo
graphs}
\author{{\small \ Fatih YILMAZ\thanks{%
e-mail adresses: fyilmaz@selcuk.edu.tr, dbozkurt@selcuk.edu.tr}, Durmus
BOZKURT} \\
{\small Selcuk University, Science Faculty Department of Mathematics, 42250 }%
\linebreak \\
{\small Campus} {\small Konya, Turkey}}
\maketitle

\begin{abstract}
In the present paper, we define two directed pseudo graph. Then, we
investigate the adjacency matrices of the defined graphs and show that the
permanents of the adjacency matrices are Jacobsthal and Jacobsthal-Lucas
numbers. We also give complex factorization formulas for the Jacobsthal
sequence.
\end{abstract}

\section{Introduction}

The \textit{Jacobsthal}\ and \textit{Jacobsthal-Lucas} sequences are defined
by the following recurrence relations, respectively: 
\begin{eqnarray*}
J_{n+2} &=&J_{n+1}+2J_{n}\text{ \ \ where }J_{0}=0,~J_{1}=1,\  \\
&& \\
j_{n+2} &=&j_{n+1}+2j_{n},\text{ \ \ where }j_{0}=2,~j_{1}=1,~\ 
\end{eqnarray*}%
for $n\geq 0.$ The first few values of the sequences are given by the
following table;%
\begin{equation*}
\begin{tabular}{c|ccccccccc}
$n$ & $1$ & $2$ & $3$ & $4$ & $5$ & $6$ & $7$ & $8$ & $9$ \\ \hline
$J_{n}$ & $1$ & $1$ & $3$ & $5$ & $11$ & $21$ & $43$ & $85$ & $171$ \\ 
$j_{n}$ & $1$ & $5$ & $7$ & $17$ & $31$ & $65$ & $127$ & $257$ & $511$%
\end{tabular}%
.
\end{equation*}

The \textit{permanent} of a matrix is similar to the determinant but all of
the signs used in the Laplace expansion of minors are positive. The
permanent of an $n$-square matrix is defined by%
\begin{equation*}
perA=\underset{\sigma \in S_{n}}{\dsum }\underset{i=1}{\overset{n}{\dprod }}%
a_{i\sigma (i)}
\end{equation*}%
where the summation extends over all permutations $\sigma $ of the symmetric
group $S_{n}$ $\cite{5}.$

Let $A=[a_{ij}]$ be an $m\times n$ matrix with row vectors $%
r_{1},r_{2},\ldots ,r_{m}.$ We call $A$ \textit{contractible} on column $k$,
if column $k$ contains exactly two non zero elements. Suppose that $A$ is
contractible on column $k$ with $a_{ik}\neq 0\neq a_{jk}$ and $i\neq j.$
Then the $(m-1)\times (n-1)$ matrix $A_{ij:k}$ obtained from $A$ replacing
row $i$ with $a_{jk}r_{i}+a_{ik}r_{j}$ and deleting row $j$ and column $k$
is called the \textit{contraction} of $A$ on column $k$ relative to rows $i$
and $j$. If $A$ is contractible on row $k$ with $a_{ki}\neq 0\neq a_{kj}$
and $i\neq j,$ then the matrix $A_{k:ij}=[A_{ij:k}^{T}]^{T}$ is called the
contraction of $A$ on row $k$ relative to columns $i$ and $j$. We know that
if $A$ is a nonnegative matrix and $B$ is a contraction of $A~\cite{2}$,
then 
\begin{equation}
perA=perB.  \label{10}
\end{equation}

A \textit{directed pseudo graph} $G=(V,E),$ with set of vertices $V(G)$ $%
=\{1,2,\ldots ,n\}$ and set of edges $E(G)=\{e_{1},e_{2},\ldots ,e_{m}\}$,
is a graph in which loops and multiple edges are allowed. A directed graph
represented with arrows on its edges, each arrow pointing towards the head
of the corresponding arc. The \textit{adjacency matrix} $A(G)=[a_{i,j}]$ is $%
n\times n$ matrix, defined by the rows and the columns of $A(G)$ are indexed
by $V(G)$, in which a$_{i,j}$ is the number of edges joinning $v_{i}$ and $%
v_{j}$ $\cite{8}$.

It is known that there are a lot of relations between determinants or
permanents of matrices and well-known number sequences. For example, in $%
\cite{6},$ the authors investigated the relations between Hessenberg
matrices and Pell and Perrin numbers. In $\cite{3},$ the authors gave a
relation between determinants of tridiagonal matrices and Lucas sequence.
Secondly, they obtained a complex factorization formula for Lucas sequence.
In $\cite{2},$ the authors considered the relationships between the sums of
Fibonacci and Lucas numbers by Hessenberg matrices.

In $\cite{4},$ the author investigated general tridiagonal matrix
determinants and permanents. Also he showed that the permanent of the
tridiagonal matrix based on $\{a_{i}\}$, $\{b_{i}\},$ $\{c_{i}\}$ is equal
to the determinant of the matrix based on $\{-a_{i}\}$, $\{b_{i}\},$ $%
\{c_{i}\}.$

In $\cite{7},~$the authors found $(0,1,-1)$ tridiagonal matrices whose
determinants and permanents are negatively subscripted Fibonacci and Lucas
numbers. Also, they give an $n\times n$ $(1,-1)$ matrix $S$, such that per$A$%
=det$(A\circ S$), where $A\circ S$ denotes Hadamard product of $A$ and $S$.
Let $S$ be a $(1,-1)$ matrix of order $n$, defined with 
\begin{equation}
S=\left[ 
\begin{array}{rrrrc}
1 & 1 & \cdots & 1 & 1 \\ 
-1 & 1 & \cdots & 1 & 1 \\ 
1 & -1 & \cdots & 1 & 1 \\ 
\vdots & \vdots & \ddots & \vdots & \vdots \\ 
1 & 1 & \cdots & -1 & 1%
\end{array}%
\right] .  \label{15}
\end{equation}

In $\cite{8},$ the authors investigate Jacobsthal numbers and obtain some
properties for the Jacobsthal numbers. They also give Cassini-like formulas
for Jacobsthal numbers as: 
\begin{equation}
J_{n+1}J_{n-1}-J_{n}^{2}=(-1)^{n}2^{n-1}  \label{12}
\end{equation}

In \cite{9}$,$ the authors investigate incomplete Jacobsthal and
Jacobsthal-Lucas numbers.

In $\cite{10},$ the author derived an explicit formula which corresponds to
the Fibonacci numbers for the number of spanning trees given below:

\FRAME{dhFw}{2.3601in}{1.2099in}{0pt}{}{}{tree1.eps}{\special{language
"Scientific Word";type "GRAPHIC";maintain-aspect-ratio TRUE;display
"USEDEF";valid_file "F";width 2.3601in;height 1.2099in;depth
0pt;original-width 11.1976in;original-height 5.7242in;cropleft "0";croptop
"1";cropright "1";cropbottom "0";filename
'çalýþmalarým/fib.blok-2k/tree1.eps';file-properties "XNPEU";}}

In $\cite{11},$ the authors consider the number of independent sets in
graphs with two elementary cycles. They described the extremal values of the
number of independent sets using Fibonacci and Lucas numbers. In $\cite{12},$
the authors give a generalization for known-sequences and then they give the
graph representations of the sequences. They generalize Fibonacci, Lucas,
Pell and Tribonacci numbers and they show that the sequences are equal to
the total number of $k$-independent sets of special graphs.

In $\cite{13},$ the authors present a combinatorial proof that the wheel $%
W_{n}$ has $L_{2n}-2$ spanning trees, L$_{n}$ is the nth Lucas number and
that the number of spanning trees of a related graph is a Fibonacci number.

In $\cite{14},$ the authors consider certain generalizations of the
well-known Fibonacci and Lucas numbers, the generalized Fibonacci and Lucas $%
p$-numbers. Then they give relationships between the generalized Fibonacci $%
p $-numbers $F_{p}(n)$, and their sums, $\overset{n}{\underset{i=1}{\dsum }}%
F_{p}(i)$, and the 1-factors of a class of bipartite graphs. Further they
determine certain matrices whose permanents generate the Lucas $p$-numbers
and their sums.

In $\cite{15},$ Lee considered $k$-Lucas and $k$-Fibonacci sequences and
investigated the relationships betwen these sequences and 1-factos of a
bipartite graph.

In the present paper, we investigate relationships between adjacency
matrices of graphs and the Jacobsthal and Jacobsthal-Lucas sequences. We
also give complex factorization formulas for the Jacobsthal numbers.

\section{Determinantal representations of the \protect\linebreak Jacobsthal
and Jacobsthal-Lucas numbers}

In this section, we consider a class of pseduo graph given in Figure 1 and
Figure 2, respectively. Then we investigate relationships between permanents
of the adjacency matrices of the graphs and the Jacobsthal and
Jacobsthal-Lucas numbers.

\FRAME{dhF}{3.3373in}{0.5535in}{0pt}{}{}{bel11.eps}{\special{language
"Scientific Word";type "GRAPHIC";maintain-aspect-ratio TRUE;display
"USEDEF";valid_file "F";width 3.3373in;height 0.5535in;depth
0pt;original-width 12.5441in;original-height 2.0574in;cropleft "0";croptop
"1";cropright "1";cropbottom "0";filename 'grafik/Bel11.eps';file-properties
"XNPEU";}}

\qquad \qquad \qquad \qquad \qquad \qquad \qquad Figure 1

Let $H_{n}=[h_{ij}]_{n\times n}$ be the adjacency matrix of the graph given
by Figure 1, in which the subdiagonal entries are $1$s, the main diagonal
entries are $1$s, except the first one which is $3$, the superdiagonal
entries are $2$s and otherwise $0.$ In other words:

\begin{equation}
H_{n}=\left[ 
\begin{array}{ccccccc}
3 & 2 &  &  &  &  &  \\ 
1 & 1 & 2 &  &  & 0 &  \\ 
& 1 & 1 & 2 &  &  &  \\ 
&  & \ddots & \ddots & \ddots &  &  \\ 
&  &  & 1 & 1 & 2 &  \\ 
& 0 &  &  & 1 & 1 & 2 \\ 
&  &  &  &  & 1 & 1%
\end{array}%
\right]  \label{30}
\end{equation}

\begin{theorem}
Let $H_{n}$ be an $n$-square matrix as in (\ref{30}), then 
\begin{equation*}
perH_{n}=perH_{n}^{(n-2)}=J_{n+2}
\end{equation*}%
where $J_{n}$ is the $n$th Jacobsthal number.
\end{theorem}

\begin{proof}
By definition of the matrix $H_{n},$ it can be contracted on column 1. Let $%
H_{n}^{(r)}$ be the $r$th contraction of $H_{n}$. If $r=1$, then 
\begin{equation*}
H_{n}^{(1)}=\left[ 
\begin{array}{cccccc}
5 & 6 &  &  &  & 0 \\ 
1 & 1 & 2 &  &  &  \\ 
& 1 & 1 & 2 &  &  \\ 
&  & \ddots & \ddots & \ddots &  \\ 
&  &  & 1 & 1 & 2 \\ 
0 &  &  &  & 1 & 1%
\end{array}%
\right] .
\end{equation*}%
Since $H_{n}^{(1)}$ also can be contracted according to the first column, we
obtain;%
\begin{equation*}
H_{n}^{(2)}=\left[ 
\begin{array}{cccccc}
11 & 10 &  &  &  & 0 \\ 
1 & 1 & 2 &  &  &  \\ 
& 1 & 1 & 2 &  &  \\ 
&  & \ddots & \ddots & \ddots &  \\ 
&  &  & 1 & 1 & 2 \\ 
0 &  &  &  & 1 & 1%
\end{array}%
\right] .
\end{equation*}%
Going with this process, we have%
\begin{equation*}
H_{n}^{(3)}=\left[ 
\begin{array}{cccccc}
21 & 22 &  &  &  & 0 \\ 
1 & 1 & 2 &  &  &  \\ 
& 1 & 1 & 2 &  &  \\ 
&  & \ddots & \ddots & \ddots &  \\ 
&  &  & 1 & 1 & 2 \\ 
0 &  &  &  & 1 & 1%
\end{array}%
\right] .
\end{equation*}%
Continuing this method, we obtain the $r$th contraction%
\begin{equation*}
H_{n}^{(r)}=\left[ 
\begin{array}{cccccc}
J_{r+2} & 2J_{r+1} &  &  &  & 0 \\ 
1 & 1 & 2 &  &  &  \\ 
& 1 & 1 & 2 &  &  \\ 
&  & \ddots & \ddots & \ddots &  \\ 
&  &  & 1 & 1 & 2 \\ 
0 &  &  &  & 1 & 1%
\end{array}%
\right]
\end{equation*}%
where $2\leq r\leq n-4.$ Hence; 
\begin{equation*}
H_{n}^{(n-3)}=\left[ 
\begin{array}{ccc}
J_{n-1} & 2J_{n-2} & 0 \\ 
1 & 1 & 2 \\ 
0 & 1 & 1%
\end{array}%
\right] \text{ }
\end{equation*}%
which, by contraction of $H_{n}^{(n-3)}$ on column 1, 
\begin{equation*}
H_{n}^{(n-2)}=\left[ 
\begin{array}{cc}
J_{n} & 2J_{n-1} \\ 
1 & 1%
\end{array}%
\right] .
\end{equation*}%
By (\ref{10}), we have $perH_{n}=perH_{n}^{(n-2)}=J_{n+2}$.
\end{proof}

Let $K_{n}=[k_{ij}]_{n\times n}$ be the adjacency matrix of the pseudo graph
given in Figure 2, with subdiagonal entries are $1$s, the main diagonal
entries are $1$s, except the second one which is $3$, the superdiagonal
entries are $2$s and otherwise $0.$ That is:

\FRAME{ftbphF}{3.2318in}{0.6806in}{0pt}{}{}{bel12.eps}{\special{language
"Scientific Word";type "GRAPHIC";maintain-aspect-ratio TRUE;display
"USEDEF";valid_file "F";width 3.2318in;height 0.6806in;depth
0pt;original-width 12.1273in;original-height 2.5313in;cropleft "0";croptop
"1";cropright "1";cropbottom "0";filename 'grafik/Bel12.eps';file-properties
"XNPEU";}}

\qquad \qquad \qquad \qquad \qquad \qquad \qquad Figure 2

and the adjacency matrix is:%
\begin{equation}
K_{n}=\left[ 
\begin{array}{ccccccc}
1 & 2 &  &  &  &  &  \\ 
1 & 3 & 2 &  &  & 0 &  \\ 
& 1 & 1 & 2 &  &  &  \\ 
&  & \ddots & \ddots & \ddots &  &  \\ 
&  &  & 1 & 1 & 2 &  \\ 
& 0 &  &  & 1 & 1 & 2 \\ 
&  &  &  &  & 1 & 1%
\end{array}%
\right]  \label{40}
\end{equation}

\begin{theorem}
Let $K_{n}$ be an $n$-square matrix as in (\ref{40}), then 
\begin{equation*}
perK_{n}=perK_{n}^{(n-2)}=j_{n}
\end{equation*}%
where $j_{n}$ is the $n$th Jacobsthal-Lucas number.
\end{theorem}

\begin{proof}
By definition of the matrix $K_{n}$, it can be contracted on column $1$.
Namely, 
\begin{equation*}
K_{n}^{(1)}=\left[ 
\begin{array}{cccccc}
5 & 2 &  &  &  & 0 \\ 
1 & 1 & 2 &  &  &  \\ 
0 & 1 & 1 & 2 &  &  \\ 
&  & \ddots & \ddots & \ddots &  \\ 
&  &  & 1 & 1 & 2 \\ 
0 &  &  &  & 1 & 1%
\end{array}%
\right] .
\end{equation*}%
$K_{n}^{(1)}$ also can be contracted on the first column, so we get;%
\begin{equation*}
K_{n}^{(2)}=\left[ 
\begin{array}{cccccc}
7 & 10 &  &  &  & 0 \\ 
1 & 1 & 2 &  &  &  \\ 
0 & 1 & 1 & 2 &  &  \\ 
&  & \ddots & \ddots & \ddots &  \\ 
&  &  & 1 & 1 & 2 \\ 
0 &  &  &  & 1 & 1%
\end{array}%
\right] .
\end{equation*}%
Continuing this process, we have%
\begin{equation*}
K_{n}^{(r)}=\left[ 
\begin{array}{cccccc}
j_{r+2} & 2j_{r+1} &  &  &  & 0 \\ 
1 & 1 & 2 &  &  &  \\ 
0 & 1 & 1 & 2 &  &  \\ 
&  & \ddots & \ddots & \ddots &  \\ 
&  &  & 1 & 1 & 2 \\ 
0 &  &  &  & 1 & 1%
\end{array}%
\right]
\end{equation*}%
for $1\leq r\leq n-4.$ Hence 
\begin{equation*}
K_{n}^{(n-3)}=\left[ 
\begin{array}{ccc}
j_{n-2} & 2j_{n-3} & 0 \\ 
1 & 1 & 2 \\ 
0 & 1 & 1%
\end{array}%
\right]
\end{equation*}%
which by contraction of $K_{n}^{(n-3)}$ on column $1$, gives%
\begin{equation*}
K_{n}^{(n-2)}=\left[ 
\begin{array}{cc}
j_{n-1} & 2j_{n-2} \\ 
1 & 1%
\end{array}%
\right] .
\end{equation*}%
By applying (\ref{10}) we have $perK_{n}=perK_{n}^{(n-2)}=j_{n},~$which is
desired.
\end{proof}

See Appendix B.

Let $S$ be a matrix as in (\ref{15}) and denote the matrices $H_{n}\circ S$
and $K_{n}\circ S$ by $A_{n}$ and $B_{n}$, respectively. Thus 
\begin{equation}
A_{n}=\left[ 
\begin{array}{cccccc}
3 & 2 &  &  &  &  \\ 
-1 & 1 & 2 &  & 0 &  \\ 
& \ddots & \ddots & \ddots &  &  \\ 
&  & -1 & 1 & 2 &  \\ 
& 0 &  & -1 & 1 & 2 \\ 
&  &  &  & -1 & 1%
\end{array}%
\right]  \label{75}
\end{equation}%
and 
\begin{equation}
B_{n}=\left[ 
\begin{array}{cccccc}
1 & 2 &  &  &  &  \\ 
-1 & 3 & 2 &  & 0 &  \\ 
& \ddots & \ddots & \ddots &  &  \\ 
&  & -1 & 1 & 2 &  \\ 
& 0 &  & -1 & 1 & 2 \\ 
&  &  &  & -1 & 1%
\end{array}%
\right] .  \label{85}
\end{equation}%
Then, we have 
\begin{equation*}
det(A_{n})=perH_{n}=J_{n+2}
\end{equation*}%
and 
\begin{equation*}
det(B_{n})=perK_{n}=j_{n}.
\end{equation*}%
In order to check the results, a Matlab and a Maple procedures are given in
Appendix $A$ and Appendix $B$ for the matrices $H_{n}$ and $K_{n}$
respectively..

\textbf{Appendix A. }We give a Matlab source code to check permanents of the
matrices given by (4)\textbf{.}

clc; clear;

$x=[$ $];$

$n=..;$

$x=eye(n);$

$\qquad x(1,1)=3;$

$\qquad x(2,2)=1;$

$for~i=1:n-1$

\qquad\ $x(i+1,i)=1;$

\qquad\ $x(i,i+1)=2;$

end

for $i=1:n-1$

$\qquad s1=[];s2=[];$

for $j=1:n+1-i$

\qquad s$1(j)=x(1,1)\ast x(2,j);$

\qquad\ s$2(j)=x(2,1)\ast x(1,j);$

end

s$1$

s$2$

st=s$1$+s$2$

if $i\symbol{126}=n-1$

\qquad\ $xy=x(3:n+1-i,2:n+1-i);$

\qquad\ $xy=[st(:,2:n+1-i);xy];$

else

xy=st$(2)$;

end

x=xy

end

x

\textbf{Appendix B. }We also give a Maple source code to check permanents of
the matrices given by (5).\textbf{\ }

\textit{restart:}

\textit{\TEXTsymbol{>} with(LinearAlgebra):}

\textit{\TEXTsymbol{>} permanent:=proc(n)}

\textit{\TEXTsymbol{>} local i,j,k,c,C;}

\textit{\TEXTsymbol{>} c:=(i,j)-\TEXTsymbol{>}piecewise(i=j+1,1,j=i+1,2,j=2
and i=2,3,i=j,1);}

\textit{\TEXTsymbol{>} C:=Matrix(n,n,c):}

\textit{\TEXTsymbol{>} for k from 0 to n-3 do}

\textit{\TEXTsymbol{>} print(k,C):}

\textit{\TEXTsymbol{>} for j from 2 to n-k do}

\textit{\TEXTsymbol{>} C[1,j]:=C[2,1]*C[1,j]+C[1,1]*C[2,j]:}

\textit{\TEXTsymbol{>} od:}

\textit{\TEXTsymbol{>} C:=DeleteRow(DeleteColumn(Matrix(n-k,n-k,C),1),2):}

\textit{\TEXTsymbol{>} od:}

\textit{\TEXTsymbol{>} print(k,eval(C)):}

\textit{\TEXTsymbol{>} end proc:with(LinearAlgebra):}

\textit{\TEXTsymbol{>} permanent( );}

\end{document}